\documentclass[12pt,draftcls,onecolumn]{IEEEtran}
\usepackage{amsmath,amssymb,array,palatino,graphicx,verbatim}

\usepackage{latexsym,url}
\usepackage{amsmath,amssymb,amsthm}
\usepackage{algorithm,algorithmic}
\usepackage{color,array,graphicx,verbatim}
\usepackage{cite} 
\usepackage{setspace}
\usepackage
  [breaklinks,bookmarks,bookmarksnumbered,bookmarksopen,bookmarksopenlevel=2]
  {hyperref}

\newcommand{\reals}{{\mbox{\bf R}}}
\newcommand{\argmax}{\mathop{\rm argmax}}
\newcommand{\vect}[1]{\boldsymbol{#1}}
\newcommand{\ie}{{\it i.e.}}
\newcommand{\eg}{{\it e.g.}}

\newtheorem{theorem}{Theorem}
\newtheorem{lemma}[theorem]{Lemma}
\newtheorem{definition}{Definition}

\usepackage{color}
\newcommand{\mycomment}[3]%
{%
\marginpar{%
  \hfil%
  \tiny{\textcolor{#2}{{\bf\textsc{#1}}}}%
  \hfil%
}%
\footnote{\textcolor{#2}{{\bf\textsc{#1}:}~~#3}}
}
\newcommand{\kakhbod}[1]%
{\mycomment{kakhbod}{blue}{#1}}
\newcommand{\jckoo}[1]%
{\mycomment{jckoo}{red}{#1}}

\begin{document}

\title{A Taxation Policy for Maximizing Social Welfare in Networks: A General Framework}
\author{\Large{Ali~Kakhbod$^1$, Joseph~Koo$^2$ and  Demosthenis~Teneketzis$^3$}\\ 
\large{${}^{1,3}$University of Michigan, Ann Arbor, ${}^2$Stanford University} \\
\normalsize{Email: ${}^{1,3}${\tt{\{akakhbod,teneket\}@umich.edu}}, ${}^2${\tt{jckoo@stanford.edu}}}}

\maketitle



\date{Preliminary Draft, April, 2010}

\section{Model}
\label{sec:ext_mkt_overview}

\subsection{Problem Formulation}
\label{subsec:ext_mkt_formulation}

Suppose there are $n$ goods which are each infinitely divisible; let
$\mathcal{L} = \{l_1,l_2,\ldots,l_n\}$ be the set of goods.  There is
only a limited amount of each good; the maximum amount of $l_j$
available is denoted by $c_{l_j}$, for $j=1,2,\ldots,n$, and is always
nonnegative.  Furthermore, we have $m$ individuals (set of individuals
denoted by $\mathcal{I}$) partitioned into $|\mathcal{G}|$ disjoint
groups; the set of groups is $\mathcal{G} =
\{g_1,g_2,\ldots,g_{|\mathcal{G}|}\}$.  Each group $g \in \mathcal{G}$
is specified by $\{\mathcal{I}^{(g)}, \mathcal{L}^{(g)}\},$ where
$\mathcal{I}^{(g)} \subseteq \mathcal{I}$ is the set of individuals in
$g$, and $\mathcal{L}^{(g)} \subseteq \mathcal{L}$ is the set of goods
requested by members of the group.  Let $\mathcal{I}_l \subseteq
\mathcal{I}$ be the set of individuals that request good $l$, and also
assume that $\mathcal{I}_l^{(g)} = \mathcal{I}_l \cap
\mathcal{I}^{(g)}$ is the set of all individuals of group $g$
requesting good $l \in \mathcal{L}^{(g)}$.  Moreover, let
$\mathcal{G}_l \subseteq \mathcal{G}$ designate the set of groups $g$
for which $l \in \mathcal{L}^{(g)}$.

Consider a specific individual $i$, where $i \in g$ and $g \in
\mathcal{G}$.  Let $\mathcal{L}_i = \{l_{i_1}, l_{i_2}, \ldots,
l_{i_{|\mathcal{L}_i|}}\} \subseteq \mathcal{L}$ be the subset of
goods which may be requested by individual $i$ (so $\mathcal{L}^{(g)}
= \bigcup_{i \in \mathcal{I}^{(g)}}\mathcal{L}_i$).  For each
individual, the subset $\mathcal{L}_i$ is known and fixed in advance.
The amount of goods actually demanded by the individual is given by
the demand vector $\vect{x}_i = (x_{i_1}, x_{i_2}, \ldots,
x_{i_{|\mathcal{L}_i|}})$, where $x_{i_k}$ is the amount of good
$l_{i_k}$ requested by $i$, for $k = 1,2,\ldots,|\mathcal{L}_i|$.  For
a given demand $\vect{x}_i$, the utility to individual $i$ is the
function $U_i : \reals^{|\mathcal{L}_i|} \to \reals$, which is
monotone\footnotemark, concave, and satisfies $U_i(\vect{0}) = 0$ and
$U_i(\vect{z}) = -\infty$ if any entry of $\vect{z}$ is negative.
Also, call $\vect{x} = (\vect{x}_1,\ldots,\vect{x}_m)$ to be the
overall demand of all individuals.
\footnotetext{Actually, we do not need $U_i(\vect{x}_i)$ to be
monotonic, as that is not necessary for our analysis.}



Suppose that we have a designer (\eg, the government) who wants to
design a mechanism to maximize the social welfare $\sum_{i=1}^m [
U_i(\vect{x}_i) - t_i(\vect{x}_i) ]$.  Here, $t_i :
\reals^{|\mathcal{L}_i|} \to \reals$ consist of taxation policies on
individuals $i$, $i=1,\ldots,m$, and is to be designed.  We consider a
scenario where the role of the designer is purely wealth
redistributionary; there is no net tax collected, so $\sum_{i=1}^m
t_i(\vect{x}_i) = 0$.  Then we can write the tax-explicit social
welfare maximization problem as the following:
\begin{equation}
\begin{array}{ll}
\mbox{maximize} & \displaystyle \sum_{i=1}^m [ U_i(\vect{x}_i)
- t_i(\vect{x}_i) ] \\
\mbox{subject to}
& \displaystyle \sum_{i=1}^m t_i(\vect{x}_i) = 0 \\
& \displaystyle \sum_{g \in \mathcal{G}} \max_{i \in
\mathcal{I}_l^{(g)}} x_{i_l} \leq c_l, \; \; \forall l \in \mathcal{L}
\end{array}
\label{eq:tax_explicit_extended_market_problem}
\mbox{,}
\end{equation}
where the optimization variables are $\vect{x}_i \in
\reals^{|\mathcal{L}_i|}$ for all $i \in \mathcal{I}$.  The resulting
taxes charged to (or accrued by) the individuals are denoted by the
vector $\vect{t} = (t_1,t_2,\ldots,t_m)$, with $t_i$ being shorthand
for $t_i(\vect{x}_i)$.
\begin{definition}
A taxation scheme $\vect{t} = (t_1, t_2, \ldots, t_m)$ is
\emph{budget-balanced} if the sum of taxes is zero; \ie, $\sum_{i=1}^m
t_i = 0$.
\end{definition}

Notice that the preceding welfare maximization problem
\eqref{eq:tax_explicit_extended_market_problem} gives the same
solution as the following social welfare maximization problem (with no
explicit taxation term):
\begin{equation}
\begin{array}{ll}
\mbox{maximize} & \displaystyle \sum_{i=1}^m U_i(\vect{x}_i) \\
\mbox{subject to}
& \displaystyle \sum_{g \in \mathcal{G}} \max_{i \in
\mathcal{I}_l^{(g)}} x_{i_l} \leq c_l, \; \; \forall l \in \mathcal{L}
\end{array}
\label{eq:extended_market_problem_with_max}
\mbox{.}
\end{equation}

We expand the constraints of problem
\eqref{eq:extended_market_problem_with_max}, in order to aid in the
decomposition.  For each good $l \in \mathcal{L}$, define the vector
$\vect{\pi}_l = (\pi_l^{(1)}, \pi_l^{(2)}, \ldots,
\pi_l^{(|\mathcal{G}|)})$ as a selection of individuals---one
individual from each group---such that every selected individual may
request good $l$.  That is, $\pi_l^{(g)}$ denotes a particular
individual such that $\pi_l^{(g)} \in \mathcal{I}_l^{(g)}$, for every
$g \in \mathcal{G}$.  If $\mathcal{I}_l^{(g)} = \emptyset$ (\ie, no
individuals in group $g$ requests good $l$), then we can ignore the
$\pi_l^{(g)}$ entry.  (We will see shortly how this is incorporated
when solving our problem.) Then for each good $l$, the set of all
possible combinations of selecting individuals (who might demand $l$)
from the groups is $\Pi_l = \{ (\pi_l^{(1)}, \pi_l^{(2)}, \ldots,
\pi_l^{(|\mathcal{G}|)}) \; | \; \pi_l^{(g)} \in \mathcal{I}_l^{(g)},
~\forall g \in \mathcal{G} \}$.  Thus $\displaystyle |\Pi_l| =
\prod_{g \in \mathcal{G}} |\mathcal{I}_l^{(g)}|$.  For later
convenience, we denote $\displaystyle P = \sum_{l \in
\mathcal{L}} |\Pi_l|$.

Equivalent to problem \eqref{eq:extended_market_problem_with_max}, we
obtain the following problem:
\begin{equation}
\begin{array}{ll}
\mbox{maximize} & \displaystyle \sum_{i=1}^m U_i(\vect{x}_i) \\
\mbox{subject to}
& \displaystyle \sum_{g \in \mathcal{G}} x_{\pi_l^{(g)}} \leq c_l,
\; \forall l \in \mathcal{L}, \; \forall \vect{\pi}_l \in \Pi_l
\end{array}
\label{eq:extended_market_problem}
\mbox{.}
\end{equation}
We call problem \eqref{eq:extended_market_problem} the primal problem.
The contribution of this paper is a simple budget-balanced taxation
scheme which is simple to implement and achieves the maximal social
welfare.

In fact, we have the following assumptions over the information
structure:
\newcounter{Acount}  
\begin{list}{(A\arabic{Acount})}{\usecounter{Acount}%
\addtolength{\leftmargin}{-\labelwidth}%
\settowidth{\labelwidth}{(A\arabic{Acount})}%
\addtolength{\leftmargin}{\labelwidth}}
\item \label{it:assume_private_utility_fcn}
The utility function of each individual will be his own private
information and need not be known by the designer.
\item \label{it:assume_designer_knowledge}
The designer does know the set of requested goods $\mathcal{L}_i$ for
each individual $i$.  Moreover, the set $\mathcal{L}_i$ is fixed.
\item \label{it:assume_individuals_pricetaker}
The individuals are price takers.
%
%
\end{list}

\section{T\^{a}tonnement Process}
\label{sec:ext_mkt_decomp}

\subsection{Dual Decomposition}
\label{subsec:ext_mkt_decomp_dual}

We consider a decomposition of the welfare maximization problem
\eqref{eq:extended_market_problem}.  From this, we will be able to
derive the taxation policy which satisfies the tax-explicit welfare
maximization problem \eqref{eq:tax_explicit_extended_market_problem}
with the stated assumptions.

Let us consider the Lagrangian of \eqref{eq:extended_market_problem},
where the Lagrange multiplier associated with the capacity constraint
$\sum_{g \in \mathcal{G}} x_{\pi_l^{(g)}} \leq c_l$ is denoted by
$p_{l, \vect{\pi}_l}$, for each good $l \in \mathcal{L}$ and each
selector $\vect{\pi}_l \in \Pi_l$.  Let $\vect{p} \in \reals^{P}$
be the vector which consists of all the Lagrange multipliers.  The
Lagrangian is
\begin{eqnarray}
L(\vect{x},\vect{p}) 
& = & \sum_{i=1}^m U_i(\vect{x}_i) + \sum_{l \in \mathcal{L}}
      \sum_{\vect{\pi}_l \in \Pi_l} p_{l, \vect{\pi}_l} \left[
      c_l - \sum_{g \in \mathcal{G}} x_{\pi_l^{(g)}} \right]  \\ 
& = & \sum_{i=1}^m U_i(\vect{x}_i) - \sum_{i=1}^m \sum_{l \in
      \mathcal{L}_i} \sum_{\vect{\pi}_l \in \Pi_l : i \in
      \vect{\pi}_l} p_{l, \vect{\pi}_l} x_{i_l} + \sum_{l \in
      \mathcal{L}} \sum_{\vect{\pi}_l \in \Pi_l} p_{l,
      \vect{\pi}_l} c_l  \\
& = & \sum_{i=1}^m \left[ U_i(\vect{x}_i) - \sum_{l \in \mathcal{L}_i}
      \left( \sum_{\vect{\pi}_l \in \Pi_l : i \in \vect{\pi}_l}
      p_{l, \vect{\pi}_l} \right) x_{i_l} \right] + \sum_{l \in
      \mathcal{L}} \left( \sum_{\vect{\pi}_l \in \Pi_l} p_{l,
      \vect{\pi}_l} \right) c_l  \\
& = & \sum_{i=1}^m \left[ U_i(\vect{x}_i) - \sum_{l \in \mathcal{L}_i}
      p_l^i x_{i_l} \right] + \sum_{l \in \mathcal{L}} p_l
      c_l
      \label{eq:extended_market_lagrangian}
\mbox{,}
\end{eqnarray}
where we let
\begin{eqnarray}
p_l & = & \sum_{\vect{\pi}_l \in \Pi_l} p_{l,
\vect{\pi}_l}  \label{eq:def_lambda_l} \\
p_l^i & = & \sum_{\vect{\pi}_l \in \Pi_l : i \in \vect{\pi}_l}
p_{l, \vect{\pi}_l}  \label{eq:def_lambda_l_i}
\mbox{.}
\end{eqnarray}
If we define $g(\vect{p}) = \max_{\vect{x}}
L(\vect{x},\vect{p})$, then the dual problem to
\eqref{eq:extended_market_problem} is
\begin{equation}
\begin{array}{ll}
\mbox{minimize} & g(\vect{p}) \\
\mbox{subject to} & \vect{p} \geq \vect{0}
\end{array}
\label{eq:extended_market_dual_problem}
\mbox{,}
\end{equation}
with variable $\vect{p}$.  If strong duality holds (which can be
checked using a constraint qualification such as Slater's
condition~\cite{boyd:convexopt}), then the solution to the dual
problem can be used to recover the solution to the primal welfare
maximization problem.

We decompose $g(\vect{p})$ so that $g(\vect{p}) =
\sum_{i=1}^m g_i(\vect{p}) + \sum_{l \in \mathcal{L}} p_l
c_l$, where
\begin{equation}
g_i(\vect{p}) = \max_{\vect{x}_i} \left[ U_i(\vect{x}_i) -
\sum_{l \in \mathcal{L}_i} p_l^i x_{i_l} \right]
\end{equation}
for each individual $i \in \mathcal{I}$.  Then each individual $i$ can
find $g_i(\vect{p})$ as the optimal value of the following
individual subproblem:
\begin{equation}
\begin{array}{ll}
\mbox{maximize} & \displaystyle U_i(\vect{x}_i) - \sum_{l \in
\mathcal{L}_i} p_l^i x_{i_l} \end{array}
\label{eq:extended_market_individual_subproblem}
\end{equation}
(for fixed $\vect{p}$ and with variable $\vect{x}_i \in
\reals^{|\mathcal{L}_i|}$).  We denote $\bar{\vect{x}}_i =
(\bar{x}_{i_1},\ldots,\bar{x}_{i_{|\mathcal{L}_i|}})$ to be the
solution to the individual subproblem for individual $i$.  We can
readily determine that $\bar{\vect{x}}_i$ will also be the solution to
\begin{equation}
\begin{array}{ll}
\displaystyle
\mbox{maximize} & \displaystyle U_i(\vect{x}_i) - \left[ \sum_{l \in
\mathcal{L}_i} p_l^i x_{i_l} - \gamma_{\vect{p}, i}
\right]
\end{array}
\mbox{,}
\label{eq:extended_market_individual_subproblem_plus_constant}
\end{equation}
(where the variable is $\vect{x}_i \in \reals^{|\mathcal{L}_i|}$), as
long as $\gamma_{\vect{p}, i}$ is constant with respect to
$\vect{x}_i$.

We can now directly solve the dual problem
\eqref{eq:extended_market_dual_problem} by solving the following
master problem (with variable $\vect{p} \in \reals^{P}$):
\begin{equation}
\begin{array}{ll}
\mbox{minimize} & \displaystyle \sum_{i=1}^m g_i(\vect{p}) +
\sum_{l \in \mathcal{L}} p_l c_l  \\
\mbox{subject to} & \displaystyle \vect{p} \geq \vect{0}
\end{array}
\label{eq:extended_market_master_problem}
\mbox{.}
\end{equation}
The mechanism designer can solve the master problem by updating
$\vect{p}$ using the projected subgradient method
\cite{bertsekas:nonlinear_prog}.  A subgradient (with respect to
$\vect{p}$) of the objective function is $\vect{s} \in
\reals^{P}$, where the elements are given by
\begin{equation}
s_{l, \vect{\pi}_l} = c_l - \sum_{i \in \vect{\pi}_l} \bar{x}_{i_l}
\label{eq:extended_market_subgrad}
\mbox{.}
\end{equation}
Thus, at each iteration, $\vect{p}$ is updated according to
\begin{equation}
p_{l, \vect{\pi}_l} := \left[ p_{l, \vect{\pi}_l} - \beta
\left[ c_l - \sum_{i \in \vect{\pi}_l} \bar{x}_{i_l} \right]
\right]^+, \; \; \forall l \in \mathcal{L}, \; \; \forall \vect{\pi}_l
\in \Pi_l
\label{eq:extended_market_master_subgrad_update}
\mbox{,}
\end{equation}
where $[z]^+$ denotes the positive part of $z$, or that $[z]^+ = \max
\{0, z\}$.

Many techniques exist for choosing the positive step size parameter
$\beta$.\footnotemark  ~Suffice to say is that for small enough
$\beta$, convergence to the optimum of the master problem
\eqref{eq:extended_market_master_problem} is guaranteed.
\footnotetext{In this work, we have assumed that the subgradient
update step sizes are chosen appropriately so that the respective
algorithms converge.  We refer the reader to
\cite{nemirovsky:problem_complexity} for a discussion of the rates of
convergence of certain step sizes, and to
\cite{bertsekas:nonlinear_prog} for additional conditions on the step
size $\beta$  which guarantee that the optimal values of the master
dual problems will be approached.  A more thorough study of the
convergence properties of subgradient methods using both constant and
non-constant step size rules (generally using diminishing step sizes)
can also be found in \cite{nedic:incremental_subgrad_methods}.  For
our algorithms, different step size rules may be helpful for speeding
up the rate of convergence.}  

The required computations are highly decentralized.  Given $p_{l,
\vect{\pi}_l}$ where $l \in \mathcal{L}_i$ and $i \in \vect{\pi}_l$,
each individual $i$ computes his own subproblem to find
$\bar{\vect{x}}_i$.  The individual receives the $p_{l, \vect{\pi}_l}$
for which $l \in \mathcal{L}_i$ and $i \in \vect{\pi}_l$, and uses
that to determine his current demand $\bar{\vect{x}}_i$ according to
\eqref{eq:extended_market_individual_subproblem}.  On the other hand,
each $(l, \vect{\pi}_l)$ (\ie, each pairing of good $l$ and selected
individuals $\vect{\pi}_l$) can compute its own price $p_{l,
\vect{\pi}_l}$, given the relevant demands $\bar{x}_{i_l}$ of
individuals $i \in \vect{\pi}_l$ for good $l$.  Very little
information needs to be exchanged between the designer and the
individuals: From the designer (or from the ``goods''), parts of
$\vect{p}$ are sent to the appropriate individuals; from the
individuals, the demands $\bar{x}_{i_l}$ are sent back to the designer
(or to the appropriate ``goods'' $\mathcal{L}_i$ and the appropriate
group selectors).  In fact, each $(l, \vect{\pi}_l)$ does not need to
be explicitly told the $\bar{x}_{i_l}$ from each individual
individually; it only needs to measure its total demand in order to
obtain $\sum_{i \in \vect{\pi}_l} \bar{x}_{i_l}$.


\subsection{Achieving Budget-Balance}
\label{subsec:ext_mkt_decomp_budget_balance}

The form of the individual subproblem
\eqref{eq:extended_market_individual_subproblem} suggests a taxation
method which would be amenable towards achieving the global optimum of
the dual problem \eqref{eq:extended_market_dual_problem}.  As
discussed, an individual $i$ solving a subproblem of the form
\eqref{eq:extended_market_individual_subproblem_plus_constant} would
be optimal at the same demand solution $\bar{\vect{x}}_i$ as that from
\eqref{eq:extended_market_individual_subproblem}.  Thus, let us
consider tax policies of the form $t_{\vect{p}, i}(\vect{x}_i) =
\sum_{l \in \mathcal{L}_i} p_l^i x_{i_l} -
\gamma_{\vect{p}, i}$.

Consider the choice of $\gamma_{\vect{p}, i} = \sum_{l \in
\mathcal{L}_i} \sum_{\vect{\pi}_l\in \Pi_l : i\in \vect{\pi}_l}
\frac{p_{l, \vect{\pi}_l}}{|\mathcal{G}_l|-1} \sum_{\substack{j \in
\vect{\pi}_l \\ j \neq i}} \bar{x}_{j_l}$.
This choice of $\gamma_{\vect{p}, i}$ is constant with respect
to the variable $\vect{x}_i$ of individual $i$'s subproblem.  The tax
policy is then
\begin{eqnarray}
t_{\vect{p}, i} (\vect{x}_i)
& = & \sum_{l \in \mathcal{L}_i} p_l^i \bar{x}_{i_l} - \sum_{l \in
      \mathcal{L}_i} \sum_{\vect{\pi}_l \in \Pi_l : i \in
      \vect{\pi}_l} \frac{p_{l, \vect{\pi}_l}}{|\mathcal{G}_l|-1}
      \sum_{\substack{j \in \vect{\pi}_l \\ j \neq i}} \bar{x}_{j_l}  \\
& = & \sum_{l \in \mathcal{L}_i} \sum_{\vect{\pi}_l \in \Pi_l : i \in
      \vect{\pi}_l}p_{l, \vect{\pi}_l} \bar{x}_{i_l} - \sum_{l \in
      \mathcal{L}_i} \sum_{\vect{\pi}_l \in \Pi_l : i \in
      \vect{\pi}_l} \frac{p_{l, \vect{\pi}_l}}{|\mathcal{G}_l|-1}
      \sum_{\substack{j \in \vect{\pi}_l \\ j \neq i}} \bar{x}_{j_l}
\mbox{.}
\label{eq:extended_market_tax_policy}
\end{eqnarray}

We assume that $|\mathcal{G}_l| > 1$ for all $l \in \mathcal{L}$.  If
$|\mathcal{G}_{\hat{l}}| = 1$ for some good $\hat{l} \in \mathcal{L}$,
then we can \textit{a priori} set the tax rate on the good to be zero,
\ie, $p_{\hat{l}, \vect{\pi}_{\hat{l}}} = 0$, so that the tax and
demand for this good will not affect the overall budget-balance.  Even
with such restrictions, our method will still produce the optimal
social welfare maximizing solution.

\begin{lemma}
\label{thm:ext_mkt_sum_taxes_zero}
The tax policy given in \eqref{eq:extended_market_tax_policy} is
budget-balanced when $\vect{x}_i = \bar{\vect{x}}_i$ for all $i \in
\mathcal{I}$.  That is, $\sum_{i=1}^m t_{\vect{p},
i}(\bar{\vect{x}}_i) = 0$.
\begin{proof}
\begin{eqnarray}
\sum_{i=1}^m t_{\vect{p}, i}(\bar{\vect{x}}_i)
& = & \sum_{l \in \mathcal{L}} \sum_{\vect{\pi}_l \in \Pi_l} \sum_{i
      \in \vect{\pi}_l}\left[p_{l, \vect{\pi}_l} \bar{x}_{i_l} -
      \frac{p_{l, \vect{\pi}_l}}{|\mathcal{G}_l|-1} \sum_{\substack{j
      \in \vect{\pi}_l \\ j \neq i}} \bar{x}_{j_l}\right]  \\
& = & \sum_{l \in \mathcal{L}} \sum_{\vect{\pi}_l \in \Pi_l} p_{l,
      \vect{\pi}_l} \left[ \sum_{i \in \vect{\pi}_l} \bar{x}_{i_l} -
      \frac{1}{|\mathcal{G}_l|-1} \sum_{i \in \vect{\pi}_l}
      \sum_{\substack{j \in \vect{\pi}_l \\ j \neq i}} \bar{x}_{j_l}
      \right] \\
& = & 0  \label{eq:extended_market_sum_taxes_zero}
\mbox{,}
\end{eqnarray}
where \eqref{eq:extended_market_sum_taxes_zero} holds because for
every good $l \in \mathcal{L}$ and every selector $\vect{\pi}_l \in
\Pi_l$, we have $\sum_{i \in \pi_l} \bar{x}_{i_l} -
\frac{1}{|\mathcal{G}_l|-1} \sum_{i \in \vect{\pi}_l}
\sum_{\substack{j \in \vect{\pi}_l \\ j \neq i}} \bar{x}_{j_l} = 0$.

For goods $\hat{l} \in \mathcal{L}$ such that $|\mathcal{G}_{\hat{l}}|
= 1$, if we initialize $p_{\hat{l}, \vect{\pi}_{\hat{l}}} = 0$
for all $\vect{\pi}_{\hat{l}} \in \Pi_{\hat{l}}$, then the subgradient
update \eqref{eq:extended_market_master_subgrad_update} will not
deviate away from $p_{\hat{l}, \vect{\pi}_{\hat{l}}} = 0$ as
long as the only demand for $\hat{l}$ is feasible.  The taxation
policy may be slightly off-balance when $p_{\hat{l},
\vect{\pi}_{\hat{l}}} > 0$ for some $\vect{\pi}_{\hat{l}} \in
\Pi_{\hat{l}}$, but the subgradient update (and resulting tax) will
eventually force the singular individual in $\mathcal{I}_{\hat{l}}$ to
return his demand to feasibility, which will also return
$p_{\hat{l}, \vect{\pi}_{\hat{l}}} = 0$ and restore
budget-balance.
\end{proof}
\end{lemma}

One cause for concern might be how individual $i$ would obtain
knowledge of the optimal solutions $\bar{\vect{x}}_j$ for the other
individuals $j$, $j \ne i$, in order to compute the ``constant'' term
in the tax policy.  This can be decreed by the designer after every
individual has indicated his demand.  Because
\begin{equation}
\argmax_{\vect{x}_i} \left[ U_i(\vect{x}_i) - \sum_{l \in
\mathcal{L}_i} p_l^i x_{i_l} \right]
= \argmax_{\vect{x}_i} \left[ U_i(\vect{x}_i) - \left[ \sum_{l \in
\mathcal{L}_i} p_l^i x_{i_l} - \gamma_{\vect{p}, i}
\right] \right]
\mbox{,}
\end{equation}
the individual could first optimize for $U_i(\vect{x}_i) - \sum_{l \in
\mathcal{L}_i} p_l^i x_l^{i}$ to find its own $\bar{\vect{x}}_i$.
This $\bar{\vect{x}}_i$ would then be sent to the mechanism designer,
who then computes the offset $\gamma_{\vect{p}, i} = \sum_{l \in
\mathcal{L}_i} \sum_{\vect{\pi}_l \in \Pi_l : i \in \vect{\pi}_l}
\frac{p_{l, \vect{\pi}_l}}{|\mathcal{G}_l|-1} \sum_{\substack{j \in
\vect{\pi}_l \\ j \neq i}} \bar{x}_{j_l} $ for every individual $i$
and using the current $\vect{p}$.  The final tax for individual $i$
can be calculated by taking the initial $\sum_{l \in \mathcal{L}_i}
p_l^i \bar{x}_{i_l}$ and then subtracting the offset term
$\gamma_{\vect{p}, i}$ which the designer tells to him.  Each
individual's tax will then be as in
\eqref{eq:extended_market_tax_policy} and the total tax from all
individuals will be zero.

This procedure will be made more clear in the t\^{a}tonnement process
in the next section.


\subsection{T\^{a}tonnement Process}
\label{subsec:ext_mkt_decomp_implementation}

The preceding decomposition can be implemented using a tax-based
approach, as shown in Algorithm~\ref{alg:extended_market_decomp}.  The
taxation policy is explicitly given, and the individual and master
problems are clearly specified.  Here, $\epsilon > 0$ is some
appropriately-chosen convergence threshold, and the norm $\|\cdot\|$
in the convergence criterion is the $\ell_2$-norm.

\begin{algorithm}
\caption{T\^{a}tonnement process for budget-balanced welfare
maximization.}
\label{alg:extended_market_decomp}
\begin{algorithmic}[1]
\STATE \label{it:ext_mkt_initialize}
Initialize $\vect{p}$ to $\vect{p}(0) = \vect{0}$.  Set $k
:= 0$.
\REPEAT
\STATE \label{it:ext_mkt_tax_policy}
Using the current $\vect{p} = \vect{p}(k)$, the designer tells
individual $i \in \mathcal{I}$ the taxation weights for demanding
particular goods; that is, the individual is told $p_l^i$ (which
equals $\sum_{\vect{\pi}_l \in \Pi_l : i \in \vect{\pi}_l} p_{l,
\vect{\pi}_l}$) for all $l \in \mathcal{L}_i$.  The linear part of the
tax policy for $i$ is then
\[
\hat{t}_{\vect{p}, i}(z) = \sum_{l \in \mathcal{L}_i}
p_l^i z_l
\mbox{.}
\]
\STATE \label{it:ext_mkt_individual_subproblem}
For individual $i$, solve
\[
\begin{array}{ll}
\mbox{maximize} & U_i(\vect{x}_i) - \hat{t}_{\vect{p},
i}(\vect{x}_i)
\end{array}
\]
for variable $\vect{x}_i \in \reals^{|\mathcal{L}_i|}$.  Note that
this is the same individual subproblem as in
\eqref{eq:extended_market_individual_subproblem}.  Set the solution as
$\bar{\vect{x}}_i$.  Send the current solution $\bar{\vect{x}}_i$ to
the designer.
\STATE \label{it:ext_mkt_subgrad_update}
The designer updates $\vect{p}$ using
\[
p_{l, \vect{\pi}_l} (k+1) = \left[ p_{l, \vect{\pi}_l} (k)
- \beta^{(k)} \left[ c_l - \sum_{i \in \vect{\pi}_l} \bar{x}_{i_l}
\right] \right]^+
\]
for each good $l \in \mathcal{L}$ and each combination of individuals
$\vect{\pi}_l \in \Pi_l$.
\STATE \label{it:ext_mkt_lambda_update}
Update $k := k + 1$.  Set $\vect{p} := \vect{p}(k)$.
\UNTIL $\|\vect{p}(k+1) - \vect{p}(k)\| < \epsilon$
\STATE \label{it:ext_mkt_charge_tax}
The designer computes $\gamma_{\vect{p}, i, \bar{\vect{x}}_{-i}} =
\sum_{l \in \mathcal{L}_i} \sum_{\vect{\pi}_l \in \Pi_l : i \in
\vect{\pi}_l} \frac{p_{l, \vect{\pi}_l}}{|\mathcal{G}_l|-1}
\sum_{\substack{j \in \vect{\pi}_l \\ j \neq i}} \bar{x}_{j_l}$ for
every individual $i \in \mathcal{I}$.  Each individual $i$ is charged
the tax
\[
t_{\vect{p}, i, \bar{\vect{x}}_{-i}}(\bar{\vect{x}}_i) =
\hat{t}_{\vect{p}, i}(\bar{\vect{x}}_i) -
\gamma_{\vect{p}, i, \bar{\vect{x}}_{-i}}
\mbox{.}
\]
\STATE \label{it:ext_mkt_final_values}
Set $\vect{x}_i^\star := \bar{\vect{x}}_i$ for all $i \in
\mathcal{I}$.  Set $\vect{p}^\star := \vect{p}$.
\end{algorithmic}
\end{algorithm}


\begin{lemma}
\label{thm:ext_mkt_feasibility_global_opt}
At iteration $k+1$, the demand allocation $\vect{x}^\star$ found from
Algorithm~\ref{alg:extended_market_decomp} is no more than
$\frac{\epsilon}{\beta^{(k)}}$-infeasible.
\begin{proof}
Consider a good $l \in \mathcal{L}$ for which the demand is infeasible
for some selection of individuals $\vect{\pi}_l$, \ie, $\sum_{i \in
\vect{\pi}_l} x_{i_l}^\star > c_l$.  From the subgradient update,
step~\ref{it:ext_mkt_subgrad_update}, we know that $0 \leq
-\beta^{(k)} \left( c_l - \sum_{i \in \vect{\pi}_l} x_{i_l}^\star
\right) \leq p_{l, \vect{\pi}_l} (k+1) - p_{l,
\vect{\pi}_l} (k)$.  Then the following inequalities hold:
\begin{eqnarray}
\sum_{i \in \vect{\pi}_l} x_{i_l}^\star - c_l
& \leq & \frac{1}{\beta^{(k)}} ( p_{l, \vect{\pi}_l} (k+1) -
p_{l, \vect{\pi}_l} (k) )  \\
& = & \frac{1}{\beta^{(k)}} | p_{l, \vect{\pi}_l} (k+1) -
p_{l, \vect{\pi}_l} (k) | \nonumber \\
& \leq & \frac{1}{\beta^{(k)}} \| \vect{p}(k+1) -
\vect{p}(k) \|  \label{eq:norm_greater_than_indiv_terms} \\
& < & \frac{1}{\beta^{(k)}} \epsilon
\mbox{.}
\end{eqnarray}
Thus, $\sum_{i \in \vect{\pi}_l} x_{i_l}^\star < c_l +
\frac{\epsilon}{\beta^{(k)}}$.
\end{proof}
\end{lemma}
The implication of the preceding lemma is that we can choose the
convergence criterion $\epsilon$ to be arbitrarily small, in order to
obtain guarantees on the feasibility of our solution.  In order to
exactly guarantee feasibility, we can also run the algorithm until the
$\vect{p}$ updates are no longer changing---at which point
$||\vect{p}(k+1) - \vect{p}(k)|| = 0$, so that $\sum_{i
\in \vect{\pi}_l} x_{i_l}^\star \leq c_l$ by
\eqref{eq:norm_greater_than_indiv_terms}.

\subsection{Convergence of T\^{a}tonnement Process}
\label{subsec:ext_mkt_decomp_convergence}

We now show that this particular decomposition and specified tax
policy converges to the solution of the welfare maximization problem.
\begin{lemma}
\label{thm:ext_mkt_lambda_convergence}
The sequence of iterates $\vect{p}(k)$ will converge to within
$\epsilon/2$ of the true optimal solution.  At this point,
Algorithm~\ref{alg:extended_market_decomp} will terminate, as the
convergence criterion $\| \vect{p}(k+1) - \vect{p}(k) \| <
\epsilon$ will have been reached.
\begin{proof}
By strong duality between the primal problem
\eqref{eq:extended_market_problem} and the dual problem
\eqref{eq:extended_market_dual_problem}, and by complementary
slackness with respect to the primal inequality constraints, we know
that the primal optimal solution $\tilde{\vect{x}}$ and dual optimal
solution $\tilde{\vect{p}}$ satisfy $ \tilde{p}_{l,
\vect{\pi}_l} = \left[ \tilde{p}_{l, \vect{\pi}_l} - \beta
\left[ c_l - \sum_{i \in \vect{\pi}_l} \tilde{x}_{i_l} \right]
\right]^+ $ for all $l \in \mathcal{L}$ and all $\vect{\pi}_l \in
\Pi_l$, for any $\beta > 0$.  We consider this to be a fixed point of
the subgradient iteration for $\vect{p}$.

>From \cite[Proposition\ 6.3.1]{bertsekas:nonlinear_prog}, we know that
if our step sizes $\beta^{(k)}$ satisfy
\begin{equation}
0 < \beta^{(k)} < \frac{2 (g(\vect{p}(k)) -
g(\vect{\tilde{p}}))}{\| \vect{s}^{(k)} \|^2}
\label{eq:ext_mkt_beta_stepsize}
\mbox{,}
\end{equation}
then each iterate $\vect{p}(k)$ will satisfy
\begin{equation}
\| \vect{p}(k+1) - \tilde{\vect{p}} \| < \|
\vect{p}(k) - \tilde{\vect{p}} \|
\mbox{,}
\end{equation}
\ie, the subgradient updates form a contractive map between the
iterate and an optimum.  This arises from the inequality
\begin{equation}
\| \vect{p}(k+1) - \vect{\tilde{p}} \|^2 \leq \|
\vect{p}(k) - \vect{\tilde{p}} \|^2 - 2 \beta^{(k)}
(g(\vect{p}(k)) - g(\tilde{\vect{p}})) + (\beta^{(k)})^2
\| \vect{s}^{(k)} \|^2
\label{eq:ext_mkt_lambda_contract}
\mbox{,}
\end{equation}
which depends on the definition of the subgradient.  Furthermore, this
inequality implies
\begin{equation}
g(\vect{p}(k+1)) - g(\vect{\tilde{p}}) \leq \frac{ \|
\vect{p}(0) - \vect{\tilde{p}} \|^2 + \sum_{i=0}^k
(\beta^{(i)})^2 \| \vect{s}^{(i)} \|^2}{2 \sum_{i=0}^k \beta^{(i)}}
\mbox{.}
\label{eq:ext_mkt_dual_opt_bound}
\end{equation}
If we choose step sizes $\beta^{(k)}$ which are square-summable but
not summable, \ie, $\sum_{k=0}^\infty \beta^{(k)} = \infty$ and
$\sum_{k=0}^\infty (\beta^{(k)})^2 < \infty$, then the dual objective
will converge to its optimum.  For example, we could choose
$\beta^{(k)} = \beta^{(0)}/k$, where the initial step size
$\beta^{(0)}$ is chosen to ensure \eqref{eq:ext_mkt_beta_stepsize} for
all iterations $k$.  Even if the step sizes are not chosen this way,
as long as the step sizes are square-summable, then we can guarantee a
bound on the difference from the optimum by using
\eqref{eq:ext_mkt_dual_opt_bound}.

The contractive map tells us that there exists some time step
$\hat{k}$ such that $\| \vect{p}(\hat{k}) -
\tilde{\vect{p}} \| < \epsilon/2$.  Moreover, $\|
\vect{p}(\hat{k}+1) - \tilde{\vect{p}} \| < \epsilon/2$,
so the value of the dual variable $\vect{p}(\hat{k}+1)$ is also
within $\epsilon/2$ of the optimal dual solution.  Then $\|
\vect{p}(\hat{k}+1) - \vect{p}(\hat{k}) \| < \epsilon$.
The convergence criterion for $\vect{p}$ has been reached, and
this occurs when the dual iterate is sufficiently close to the optimal
dual solution.
\end{proof}
\end{lemma}

\begin{theorem}
\label{thm:ext_mkt_convergence_global_opt}
Assuming that strong duality holds,
Algorithm~\ref{alg:extended_market_decomp} converges to the global
optimum of the welfare maximization problem
\eqref{eq:extended_market_problem}.
\begin{proof}
The subgradient update in step~\ref{it:ext_mkt_subgrad_update} will
converge to the optimal solution of the dual problem
\eqref{eq:extended_market_dual_problem}, which is also the optimal
solution of the master dual problem
\eqref{eq:extended_market_master_problem}.  We know that at
convergence, the solutions $\vect{x}_i^\star$, which are the
maximizers from step~\ref{it:ext_mkt_individual_subproblem} when
$\vect{p} = \vect{p}^\star$, are the same as the
maximizers for the subproblems
\eqref{eq:extended_market_individual_subproblem} (for every $i \in
\mathcal{I}$).  By strong duality, the dual value at the solution to
\eqref{eq:extended_market_dual_problem} is the same as the primal
optimal value for \eqref{eq:extended_market_problem}.  Because the
objective function for each subproblem
\eqref{eq:extended_market_individual_subproblem} is strictly concave,
the optimal solution for each subproblem is unique, and so
$\vect{x}_i^\star$ for all $i \in \mathcal{I}$ is the solution to the
primal problem \eqref{eq:extended_market_problem}.  Thus the algorithm
gives the demand allocation which finds the maximum social welfare.
\end{proof}
\end{theorem}

\subsection{Alternative Taxation Policy for Achieving Budget-Balance (at
Equilibrium)}
\label{subsec:decomp_alt_tax_policy}
At convergence of Algorithm~\ref{alg:extended_market_decomp} (when
each individual $i$ demands an allocation of $\vect{x}_i^\star$), if
we instead use an alternative tax policy of
\begin{equation}
\tau_{\vect{p}^\star, i}(z) = \sum_{l \in \mathcal{L}_i}
(p_l^i)^\star \left[ z_l - \frac{1}{|\mathcal{G}_l|} c_l \right]  
\mbox{,}
\end{equation}
then this tax policy will be budget-balanced.  Here,
$\vect{p}^\star$ is the optimal dual solution given at algorithm
convergence.  Budget-balance can be shown by computing the sum of
taxes:
\begin{eqnarray}
\sum_{i=1}^m \tau_{\vect{p}^\star, i}(\vect{x}_i^\star)
& = & \sum_{i=1}^m \sum_{l \in \mathcal{L}_i} (p_l^i)^\star
      \left[ x_{i_l}^\star - \frac{1}{|\mathcal{G}_l|} c_l \right] \\
& = & \sum_{l \in \mathcal{L}} \sum_{i: l \in \mathcal{L}_i}
      (p_l^i)^\star x_{i_l}^\star - \sum_{l \in \mathcal{L}}
      \sum_{i: l \in \mathcal{L}_i} \frac{1}{|\mathcal{G}_l|}
      (p_l^i)^\star c_l  \\
& = & \sum_{l \in \mathcal{L}} \sum_{i: l \in \mathcal{L}_i}
      \sum_{\vect{\pi}_l \in \Pi_l : i \in \vect{\pi}_l} p_{l,
      \vect{\pi}_l}^\star x_{i_l}^\star - \sum_{l \in \mathcal{L}}
      \sum_{i: l \in \mathcal{L}_i} \sum_{\vect{\pi}_l \in \Pi_l : i
      \in \vect{\pi}_l} \frac{1}{|\mathcal{G}_l|} p_{l,
      \vect{\pi}_l}^\star c_l  \\
& = & \sum_{l \in \mathcal{L}} \sum_{\vect{\pi}_l \in \Pi_l} \sum_{i
      \in \vect{\pi}_l} p_{l, \vect{\pi}_l}^\star x_{i_l}^\star
      - \sum_{l \in \mathcal{L}} \sum_{\vect{\pi}_l \in \Pi_l} \sum_{i
      \in \vect{\pi}_l} \frac{1}{|\mathcal{G}_l|} p_{l,
      \vect{\pi}_l}^\star c_l  \\
& = & \sum_{l \in \mathcal{L}} \sum_{\vect{\pi}_l \in \Pi_l}
      p_{l, \vect{\pi}_l}^\star \sum_{i \in \vect{\pi}_l}
      x_{i_l}^\star - \sum_{l \in \mathcal{L}} \sum_{\vect{\pi}_l \in
      \Pi_l} p_{l, \vect{\pi}_l}^\star c_l \left( \sum_{i \in
      \vect{\pi}_l} \frac{1}{|\mathcal{G}_l|} \right)  \\
& = & \sum_{l \in \mathcal{L}} \sum_{\vect{\pi}_l \in \Pi_l}
      p_{l, \vect{\pi}_l}^\star \left[ \sum_{i \in \vect{\pi}_l}
      x_{i_l}^\star - c_l \right]
      \label{eq:extended_market_sum_taxes}
\mbox{,}
\end{eqnarray}
where \eqref{eq:extended_market_sum_taxes} holds because $\sum_{i \in
\vect{\pi}_l} \frac{1}{|\mathcal{G}_l|} = 1$ for all $l \in
\mathcal{L}$.  By Theorem~\ref{thm:ext_mkt_convergence_global_opt},
each ${\vect{x}}_i^\star$ is the demand allocation which globally
maximizes the social welfare.  Because the difference $\sum_{i \in
\vect{\pi}_l} x_{i_l} - c_l$ is the constraint associated with the
Lagrange multiplier $p_{l, \vect{\pi}_l}$, by complementary
slackness~\cite{boyd:convexopt}, the product $p_{l,
\vect{\pi}_l}^\star \left[ \sum_{i \in \vect{\pi}_l} x_{i_l}^\star -
c_l \right] = 0$ for every $l \in \mathcal{L}$ and $\vect{\pi}_l \in
\Pi_l$.  Thus, $\sum_{l \in \mathcal{L}} \sum_{\vect{\pi}_l \in \Pi_l}
p_{l, \vect{\pi}_l}^\star \left[ \sum_{i \in \vect{\pi}_l}
x_{i_l}^\star - c_l \right] = 0$, and the tax policy is
budget-balanced.

When Algorithm~\ref{alg:extended_market_decomp} converges, we know
that $\vect{p}^\star$ satisfies the complementary slackness
conditions.  This tells us that if a particular good $\hat{l} \in
\mathcal{L}$ under combination $\vect{\pi}_{\hat{l}}$ is not fully
demanded, \ie, when $\sum_{i \in \vect{\pi}_{\hat{l}}}
x_{i_{\hat{l}}}^\star < c_{\hat{l}}$, then from complementary
slackness we know that $p_{\hat{l}, \vect{\pi}_{\hat{l}}}^\star
= 0$.  This means that any individual $i \in \vect{\pi}_{\hat{l}}$,
where $\hat{l} \in \mathcal{L}_i$, could increase his demand
$x_{i_{\hat{l}}}$ without any taxation penalty with regards to the
particular combination $\vect{\pi}_{\hat{l}}$; however, he will not do
so as that would decrease his own utility (recall that every
individual is already at his optimal point $\vect{x}_i^\star$).  In
fact, if every utility function were strictly increasing, then all of
the maximum good demand constraints would be satisfied with equality.
This is because every individual would always want to increase his
demands---thereby increasing his utility as long as that good has no
tax penalty---until goods can no longer be provided to him.

We have discussed that imposing a tax policy of $\tau_{\vect{p},
i}(\vect{z}) = \sum_{l \in \mathcal{L}_i} p_l^i \left[ z_l -
\frac{1}{|\mathcal{G}_l|} c_l \right]$ for every individual $i$ will
lead to a solution which is budget-balanced at optimality.  However,
any tax policy of the form $\tau_{\vect{p}, i}(\vect{z}) =
\sum_{l \in \mathcal{L}_i} p_l^i z_l - \gamma_{\vect{p},
i}$, where $\gamma_{\vect{p}, i}$ is the constant offset term in
the policy, is acceptable for budget balance---as long as
$\sum_{i=1}^m \gamma_{\vect{p}, i} = \sum_{l \in \mathcal{L}}
p_l c_l$ (so that the complementary slackness budget-balance
argument still holds).  Each individual may be given a different
constant offset for its required tax, but the algorithm will still
converge to the same solution $\vect{x}^\star$ since the offsets do
not change the demand solutions of the individual subproblems
\eqref{eq:extended_market_individual_subproblem}.  It may be useful to
consider other forms for $\gamma_{\vect{p}, i}$ to satisfy some
other desired property (for example, some notion of fairness).  For
example, one approach would be to consider constant offsets of the
form $\gamma_{\vect{p}, i} = \theta_i \sum_{l \in \mathcal{L}}
p_l c_l$, where $\sum_{i=1}^m \theta_i = 1$.

\subsection{Discussion}
\label{subsec:discussion}

\begin{definition}
Given a demand $\vect{x}$, for any good $l \in \mathcal{L}$ and group
selectors $\vect{\pi}_l, \vect{\pi}_l' \in \Pi_l$, we call
\begin{eqnarray}
\vect{\pi}_l \succ \vect{\pi}_l' \quad \mbox{if} \quad \sum_{i \in
\vect{\pi}_l} x_{i_l} > \sum_{i \in \vect{\pi}_l'} x_{i_l}
\mbox{.}
\end{eqnarray}  
Consequently, define the set $\Pi_l^\mathrm{max}$ as
\begin{eqnarray}
\Pi_l^\mathrm{max} = \{\vect{\pi}_l \in \Pi_l \; | \; \vect{\pi}_l
\succ \vect{\pi}_l' ~\mbox{for all} ~\vect{\pi}_l' \in \Pi_l
~\mbox{such that} ~\vect{\pi}_l' \neq \vect{\pi}_l \}
\mbox{.}
\end{eqnarray}
For a particular good $l$, an element of $\Pi_l^\mathrm{max}$ is
denoted by $\vect{\pi}_l^\mathrm{max}$.
\end{definition}

From the t\^{a}tonnement process specified in
section~\ref{subsec:ext_mkt_decomp_implementation}, the following
observations can be made \emph{at equilibrium}:
\newcounter{Ocount}  
\begin{list}{(O\arabic{Ocount})}{\usecounter{Ocount}%
\addtolength{\leftmargin}{-\labelwidth}%
\settowidth{\labelwidth}{(O\arabic{Ocount})}%
\addtolength{\leftmargin}{\labelwidth}}
\item \label{it:dual_vars_for_max_selector}
Suppose the utility functions are monotone.  Then for any $l \in
\mathcal{L}$ and $\vect{\pi}_l^\mathrm{max} \in \Pi_{l}^\mathrm{max}$,
we have $\sum_{i \in \vect{\pi}_l^\mathrm{max}} x_i^\star - c_l = 0$.
Thus for all $\vect{\pi}_l^\mathrm{max} \in \Pi_l^\mathrm{max}$ and
$\vect{\pi}_l \in \Pi_l \setminus \Pi_l^\mathrm{max}$, we obtain the
following properties:
\begin{itemize}
\item $\sum_{i \in \vect{\pi}_l} x_{i_l}^\star < c_l$.
\item $p_{l, \vect{\pi}_l}^\star = 0$.
\item $p_l^\star = \sum_{\vect{\pi}_l^\mathrm{max} \in
\Pi_l^\mathrm{max}} p_{l, \vect{\pi}_l^\mathrm{max}}^\star$.
\item For all $j \in \mathcal{I}_l$, if $j \not\in
\vect{\pi}_l^\mathrm{max}$ for all $\vect{\pi}_l^\mathrm{max}$, then
$(p_l^j)^\star = 0$.
\item For all $j \in \mathcal{I}_l$, if there exists $\pi_l^{\max}\in \Pi_l^{\max}$ where $j \in \pi_l^{\max}$, then $(p_l^j)^\star=\sum_{\vect{\pi}_l^\mathrm{max} \in
\Pi_l^\mathrm{max}: j \in \vect{\pi}_l^\mathrm{max}} p_{l, \vect{\pi}_l^\mathrm{max}}^\star$.
\end{itemize}
\item \label{it:single_max_selector}
Consider $l \in \mathcal{L}$ and suppose that
$|\Pi_l^\mathrm{max}|=1$, \ie, $\Pi_l^\mathrm{max} =
\{\vect{\pi}_l^\mathrm{max}\}$, then the following hold:
\begin{itemize}
\item $p_l^\star = (p_l^i)^\star = (p_l^j)^\star$
for any $i,j \in \vect{\pi}_l^\mathrm{max}$.
\item For good $l$, the problem becomes a market problem for the
individuals who are in $\vect{\pi}_l^\mathrm{max}$.
\item For any $j \in \mathcal{I}_l \setminus
\{\vect{\pi}_l^\mathrm{max}\}$, $(p_l^j)^\star = 0$.
\end{itemize}
\end{list}

\section{Conclusion}
\label{sec:conclusion}
We have presented a simple t\^{a}tonnement process based on a
decomposition method which is simple to implement and achieves the
maximal social welfare, under the assumption that the utility function
of each [price-taking] individual will be his own private information
and need not be known by the designer.  At each iteration, very little
information needs to be exchanged among the individuals in order to
achieve the optimal allocation.  Furthermore, the given
t\^{a}tonnement process is always balanced at equilibrium and off equilibrium.

\bibliographystyle{IEEEtran} 
\bibliography{IEEEabrv,references}

\end{document}